\documentclass[12pt] {article}
\usepackage{
        amssymb,
        amsthm,
        amsmath
}

\newcommand {\bbR}    {\mathbb{R}}
\newcommand {\F} {F_{z_1,z_2}}
\newcommand {\arccosh}{\mathop{\mathrm{arccosh}}}

\newtheorem{theorem}{Theorem}[section]

\newtheorem{corollary}[theorem]{Corollary}

\newtheorem{remark}{Remark}[section]

\title{Dimension Reduction for the Hyperbolic Space}
\author{Itai Benjamini\thanks{Microsoft Research and The Weizmann Institute, itai.benjamini@weizmann.ac.il.}
 \and%
 Yury Makarychev\thanks{Microsoft Research, One Microsoft Way, Redmond,WA 98052. yurym@microsoft.com.}
 }

\begin{document}
\maketitle
\begin{abstract}
A dimension reduction for the hyperbolic space is established.
When points are far apart an embedding with bounded distortion
into $ H^2$ is achieved.
\end{abstract}

\section{Introduction}
Dimension reduction algorithms for Euclidean spaces have numerous 
algorithmic applications. They help to significantly reduce
the space required for storing multidimensional data, 
and thus to improve performance of many algorithms.
In this paper, we present a dimension reduction algorithm 
for the hyperbolic space. Our results show that many existing
algorithms for Euclidean spaces that rely on dimension reduction
can be also applied to hyperbolic spaces. 
We refer the reader to a paper of Ailon and Chazelle~\cite{AC} 
for background on dimension reduction algorithms and some 
of their applications. We also refer the reader to a paper of 
Krauthgamer and Lee~\cite{KL}, which studies combinatorial
algorithms for hyperbolic spaces.

For background on hyperbolic geometry see e.g.~\cite{CFKP}. 
For a recent study of convexity and high dimensional 
hyperbolic spaces see~\cite{R}.

\subsection{Our Results}
In this paper, we consider the Poincar\'e half-space model of the
hyperbolic space $ H^n$. Recall that every point is represented as
a pair $(z, x)$, $z\in \bbR^+$, $x \in \bbR^{n-1}$ in this model.
The distance between two points $p_1=(z_1,x_1)$ and $p_2=(z_2,
x_2)$ is defined by
$$d(p_1,p_2) = \arccosh \left(1 + \frac{\|x_1 - x_2\|^2 +
  (z_1 - z_2)^2}{2z_1z_2}\right).$$
For brevity, we define $F(r,z_1, z_2)$ as follows:
$$\F(r) = \arccosh \left(1 + \frac{r^2 + (z_1 - z_2)^2}{2z_1z_2}\right).$$
Then
$$d(p_1,p_2) = \F (\|x_1 -x_2\|).$$

Suppose we are given an $n$-point subset $S$ of the hyperbolic space.
Let $T$ be its projection on
$\bbR^{n-1}$:
$$T = \{x: (z,x)\in S\}.$$
By the Johnson--Lindenstrauss lemma \cite{JL},
there exists an embedding of $T$ into $O((\log n)/\varepsilon^2)$ dimensional
Euclidean space such that for every
$x_1, x_2\in T$
$$\|x_1 - x_2\| \leq
\|f(x_1)- f(x_2)\| \leq (1 + \varepsilon) \|x_1 - x_2\|.$$

\begin{theorem}[Dimension Reduction for $ H^n$] 
\label{thm:one}
Consider the map $g: H^n\to  H^{O(\log n)}$ defined by
$$g(p) \equiv g((z,x)) = (z, f(x)).$$
 Then for every two points $p_1$
and $p_2$ at (hyperbolic) distance $\Delta$, we have
$$\Delta \leq d(g(p_1), g(p_2)) \leq
\left(1 + \frac{3\varepsilon}{1+\Delta}\right) \Delta.$$
\end{theorem}
\begin{remark}
Since we reduce the hyperbolic case to the Euclidean case,
the dimension reduction embedding for $H^n$
can be computed very efficiently 
using \textit{the Fast Johnson--Lindenstrauss Transform}
of Ailon and Chazelle~\cite{AC}.
\end{remark}
The following corollary follows from a result of Bonk and
Schramm~\cite{BS}.
\begin{corollary}
Let X be a Gromov hyperbolic geodesic metric space with bounded
growth at some scale. Then there exist constants $\lambda_X$ and
$C_X$ such that every $n$-point subset $S$ of $X$ roughly
quasi-similar embeds into a $O((\log n)/\varepsilon^2)$
dimensional hyperbolic space. That is, there exists a map
$\varphi:S\to  H^{O((\log n)/\varepsilon^2)}$ such that for every
$x,y \in S$
$$\lambda_X d(x,y) - C_X \leq d(\varphi(x),\varphi(y)) \leq (1+\varepsilon)
\lambda_X d(x,y) + C_X.$$
\end{corollary}

For far apart points we prove the following theorem.

\begin{theorem}[Embedding into Hyperbolic Plane]
\label{thm:two}
Let $S$ be an $n$ point subset of $ H^n$. Assume that the distance
between every two points in $S$ is at least $\frac{\ln
(12n)}{\varepsilon}$ then there exists an embedding of $S$ into
the hyperbolic plane $ H^2$ with distortion at most
$1+\varepsilon$.
%Moreover, $\|g\|_{Lip} \leq 1 + \varepsilon$,
%$\|g\|_{Lip} \geq 1$.
\end{theorem}

\section{Proofs}
We start with the proof of the first theorem followed by the proof
for the second.

\begin{proof}[Proof of Theorem~\ref{thm:one}]
First, since $\F$ is an increasing function, we have
$$d(g(p_1), g(p_2)) = \F(\|f(x_1) - f(x_2)\|) \geq \F(\|x_1 - x_2\|) = \Delta.$$
On the other hand, by the mean value theorem,
\begin{align}
d(g(p_1), g(p_2)) & \leq \F((1+\varepsilon) \|x_1 - x_2\|) \\
&= \F(\|x_1 - x_2\|) + \frac{d\F(\hat r)}{dr} \cdot \varepsilon \|x_1 - x_2\|,
\label{eqn:ubound}
\end{align}
where $\hat r \in (\|x_1-x_2\|, (1+\varepsilon) \|x_1 - x_2\|)$.
Let us now bound the derivative of $\F$.
\begin{align*}
\frac{d\F(\hat r)}{dr} &=
{\frac{2\hat r}{2z_1 z_2}} \cdot \left.\frac{1}{\sqrt{t-1}\sqrt{t+1}}
\right|_{t=1+\frac{\hat r^2 + (z_1 - z_2)^2}{2z_1z_2}}
\\
&=\frac{2\hat r}{\sqrt{\hat r^2 + |z_1-z_2|^2}\sqrt{\hat r^2 +
|z_1-z_2|^2 + 4z_1z_2}}\\
&\leq \frac{2}{\sqrt{\hat r^2 + |z_1-z_2|^2 + 4z_1z_2}}\\
&\leq \frac{2}{\sqrt{\|x_1-x_2\|^2 + |z_1-z_2|^2 + 4z_1z_2}}.
\end{align*}
Here, we used that $(\arccosh t)' = 1/\sqrt{(t - 1)(t+1)}$.
From the identity
$$ \frac{\|x_1-x_2\|^2 + |z_1-z_2|^2}{2z_1z_2} = \cosh \Delta - 1,$$
and the bound for $\frac{d\F(\hat r)}{dr}$ we get an estimate for the
additive term in (\ref{eqn:ubound})
\begin{align*}
\frac{d\F(\hat r)}{dr} \cdot \varepsilon\|x_1 - x_2\| &\leq
\frac{2\|x_1 - x_2\|\, \varepsilon}{\sqrt{2 z_1 z_2 (\cosh \Delta + 1)}}\\
&\leq
2  \varepsilon \sqrt{\frac{\|x_1 - x_2\|^2 + |z_1-z_2|^2}{2 z_1 z_2 }}
\cdot \frac{1}{\sqrt{\cosh \Delta + 1}}\\
&= 2 \varepsilon \sqrt{\frac{\cosh \Delta - 1}{\cosh \Delta + 1}} =
2\varepsilon \tanh\frac{\Delta}{2}.
\end{align*}
It is easy to see that
$$\tanh t \leq \frac{3t}{1+ 2t}$$
for $t>0$. Therefore, the additive term in (\ref{eqn:ubound})  is at most
$$\frac{3\varepsilon}{1 + \Delta} \Delta.$$
This concludes the proof.
\end{proof}
\medskip

\begin{proof}[Proof of Theorem~\ref{thm:two}]
Define $T = \{x: (z,x) \in S\}$. By a theorem of
Matou\v{s}ek~\cite{M}, there exists an embedding $f:T\to \bbR$ of
$\bbR^{n-1}$ into $\bbR$ with distortion at most $12n$. We assume
that $f$ is non-contracting and $\|f\|_{Lip} \leq 12 n$. Consider
the embedding $g:S\to  H^2$ defined by
$$g((z,x)) = (z, f(x)).$$
Clearly, $g$ is non-contracting. Now we upper bound the Lipschitz norm
of $g$. Pick two points $p_1 =(z_1, x_1)$ and $p_2=(z_2, x_2)$ at
distance $\Delta$ in $S$. Let $r = \|x_1 - x_2\|$.
\begin{align*}
d(g(p_1), g(p_2)) &= \F(\|f(x_1)-f(x_2)\|) \leq \F(12nr)\\
&\leq \arccosh\left(1 + 12n \frac{r^2 + |z_1-z_2|^2}{2z_1 z_2}\right).
\end{align*}
Since
$$\frac{r^2 + |z_1-z_2|^2}{2z_1 z_2} = \cosh \Delta - 1,$$
we have
$$d(g(p_1), g(p_2)) \leq \arccosh(12n\cosh \Delta - (12n-1)).$$
Observe that
\begin{align*}
\cosh t &= \frac{e^t + e^{-t}}{2}\leq \frac{e^t}{2} + \frac{1}{2}\qquad (\text{for } t>0);\\
\arccosh t &= \ln( t + \sqrt{t^2-1}) \leq \ln (2t)\qquad
(\text{for } t>1).
\end{align*}
Therefore,
\begin{align*}
d(g(p_1), g(p_2))  &\leq \arccosh(12n\cosh \Delta - (12n-1)) \leq
\ln(2\cdot 12n \frac{e^\Delta}{2})\\
&= \Delta + \ln (12n) \leq (1  +\varepsilon) \Delta.
\end{align*}

This concludes the proof.
\end{proof}

\noindent
{\bf Acknowledgements:} Thanks to James Lee, Assaf Naor and Oded Schramm for
useful discussions.

\end{document}